\documentclass[12pt]{article}
\usepackage{amssymb,amsmath}

\newcommand{\w}{\wedge}
\newcommand{\n}{\notag}
\newcommand{\im}{\mbox{i}}
\newcommand{\hook}{\hookrightarrow}
\newcommand{\tb}{\textbf}
\newcommand{\vsp}{\vspace}
\newcommand{\sq}{$\; \, \square$}
\newcommand{\Remark}{\emph{Remark.}$\;$}

\newcommand{\B}{ \textbf{B} }

\newcommand{\C}{ \mathbb{C} }

\newcommand{\ap}{\alpha}

\begin{document}

\sloppy

\begin{center}

\begin{center}
\Large{\tb{A gap rigidity for proper holomorphic maps\\
from $\, \B^{n+1}$ to $\, \B^{3n-1}$}}
\end{center}

\vspace{2pc}

 Sung Ho Wang\\
 Department of Mathematics\\
 Kias \\
 Seoul, Corea 130-722 \\
\texttt{shw@kias.re.kr}
\end{center}

\vspace{3pc}

\noindent \textbf{Abstract}$\;\,$
Let $\, \B^{n+1} \subset \C^{n+1}$ be the unit ball in a complex
Euclidean space, and let $\, \Sigma^n = \partial \B^{n+1} = S^{2n+1}$.
Let $\, f: \Sigma^n \hook \Sigma^{N}$ be a local CR immersion.
If $\, N-n<2n-1$, the asymptotic vectors of the second fundamental form
of $\, f$ at each point
form a subspace of the holomorphic tangent space of
$\, \Sigma^n$ of codimension at most 1.
We exploit the successive derivatives of this relation and show that
a linearly full local CR immersion
$\, f: \Sigma^n \hook \Sigma^{N}$, $\, N \leq 3n-2$,
can only occur when $\, N = n, \, 2n$, or $\, 2n+1$.
Together with the recent classification of the rational proper holomorphic maps
from $\, \B^{n+1}$ to $\, \B^{2n+2}$ by Hamada,
this gives a classification of the rational proper holomorphic maps
from $\, \B^{n+1}$ to $\, \B^{3n-1}$ for $\, n \geq 3$.

\vspace{1pc}


\noindent \textbf{MS classification:} 32H02

\thispagestyle{empty}

\newpage
\setcounter{page}{1}

\begin{center}
\Large{\tb{A gap rigidity for proper holomorphic maps\\
from $\, \B^{n+1}$ to $\, \B^{3n-1}$}}
\end{center}

\begin{center}
Sung Ho Wang
\end{center}

\vsp{2pc}
\begin{center}
\textbf{\large{Introduction}}
\end{center}
The purpose of this paper is to prove the following gap rigidity
for the proper holomorphic maps between the unit balls
in complex Euclidean spaces.
Let $\,z= ( z^0, \, z^i)$, $\, 1 \leq i \leq n$, be the coordinates
of $\, \C^{n+1}$.
Let $\, \B^{n+1} \subset \C^{n+1}$ denote the unit ball.

\vsp{1pc}
\tb{Theorem}
\emph{
Let $\, F: \B^{n+1} \to \B^{3n-1}$ be a proper holomorphic map
that is $\, C^3$ up to the boundary, $\, n \geq 3$.
Then, up to automorphisms of the unit balls,
$\, F$ is equivalent to one of the following three polynomial maps.}

\emph{
\tb{A}. Linear embedding
$\, F: \B^{n+1} \to \B^{n+1} \subset \B^{3n-1}$}.

\emph{
\tb{B}. Whitney map
$\, F: \B^{n+1} \to \B^{2n+1} \subset \B^{3n-1}$ defined by}

\qquad \qquad \qquad
$\, F(z) = ( z^i,\, z^i  z^0, \, (z^0)^2)$.

\emph{
\tb{C}. $\, F: \B^{n+1} \to \B^{2n+2} \subset \B^{3n-1}$ is defined
for some $\, \phi, \, 0 < \phi < \frac{\pi}{2},$ by}

\qquad \qquad \qquad
$\, F(z) =
( z^i,\, \cos(\phi) \,z^0, \, \sin(\phi) \,z^i  z^0, \, \sin(\phi) \,(z^0)^2)$.
\vsp{1pc}

Proper holomorphic map between the unit balls is a subject
with a long and fruitful history that goes back to Poincre, and Alexander,
[Fo][DA2] for general references.
More recently, Huang and Ji showed that every rational proper holomorphic maps
from $\, \B^{n+1}$  to $\, \B^{2n+1}$, $\, n \geq 2$, is equivalent to
either the linear embedding or the Whitney map [HJ].
Hamada showed that every rational proper holomorphic maps
from $\, \B^{n+1}$  to $\, \B^{2n+2}$, $\, n \geq 3$, is equivalent to
either one of the maps \tb{A}, \tb{B}, or \tb{C} [Ha].
The one parameter family of maps \tb{C} was introduced by D'Angelo [DA1].

Let $\, \Sigma^n = \partial \B^{n+1} = S^{2n+1}$.
Since a proper holomorphic map $\, F: \B^{n+1} \to  \B^{3n-1}$ that is
$\, C^3$ up to the boundary is rational [HJX, Corollary 1.4],
in order to prove \tb{Theorem} it suffices to show that
any local CR immersion
$\, f: \Sigma^n \hook \Sigma^{3n-2}$ lies in
$\, f: \Sigma^n \hook \Sigma^{2n+1}\subset \Sigma^{3n-2}$ for $\, n \geq 3$,
thereby reducing the problem to the case treated by Hamada.
In fact, our analysis shows that a linearly full local CR immersion
$\, f: \Sigma^n \hook \Sigma^{N}$, $\, N \leq 3n-2$,
can only occur when $\, N = n, \, 2n$, or $\, 2n+1$.

For a local CR immersion $\, f:\Sigma^n \hook \Sigma^{N}$,
the asymptotic vectors of the second fundamental form of $\, f$ at each point
form a subspace of the holomorphic tangent space of $\, \Sigma^n$.
\emph{Rank} of the second fundamental form at a point is defined
as the codimension of the asymptotic subspace.
When $\, N-n < 2n-1$, an algebraic analysis of CR Gau\ss \, equation
shows that the rank is at most 1 [Iw].
The idea is then to explore the consequence of the successive derivatives
of this relation to obtain the desired gap theorem.
In this regard, a computations suggests that
this type of \emph{gap phenomena} may persist for local CR immersions
$\, f: \Sigma^n \to \Sigma^{\mu(n+1)-1}$ with second fundamental form of rank 1
for $\, \mu \leq n-1$.
We make a relevant remark at the end of section 2.

Huang introduced the notion of geometric rank of a CR map between
spheres [Hu]. A short computation shows that if a CR immersion has
second fundamental form of rank 1, then it has geometric rank 1.
We suspect that geometric rank is bounded above by
the rank of the second fundamental form at a generic point.

One of our initial motivation was whether the complete system for
CR maps in [Han] is subject to any geometric interpretation.

\vsp{1pc}
\section{CR submanifold}
We set up the basic structure equations
for CR submanifolds in spheres.
For general reference in CR geometry, [ChM][DA2][EHZ][Fo].

Let $\, \mathbb{C}^{N+1,1}$ be the complex vector space with coordinates
$\, z = (\, z^0, \, z^A, \, z^{N+1} \,)$, $\, 1 \, \leq  A \, \leq  N$,
and a Hermitian scalar product
\begin{equation}
\langle \, z, \, \bar{z} \, \rangle = z^A \, \bar{z}^A +
\mbox{i} \, (z^0 \, \bar{z}^{N+1} - z^{N+1} \, \bar{z}^0 ).\n
\end{equation}
Let $\, \Sigma^N$ be the set of equivalence classes up to scale of
null vectors with respect to this product. Let SU$(N+1,1)$ be the
group of unimodular linear transformations that leave the form $\,
\langle \, z, \, \bar{z} \rangle$ invariant. Then SU$(N+1,1)$ acts
transitively on $\, \Sigma^N$, and
\begin{equation}
p: \, \mbox{SU}(N+1,1) \to \Sigma^N = \, \mbox{SU}(N+1,1)/P \n
\end{equation}
for an appropriate  subgroup $\, P$ [ChM].

Explicitly, consider an element $\, Z=(\, Z_0, \, Z_A, \, Z_{N+1}\,)
\in \mbox{SU}(N+1,1)$ as an ordered set of $\, (N+2)$-column vectors
in $\, \mathbb{C}^{N+1,1}$  such that det$(Z)=1$, and that
\begin{equation}\label{prod}
\langle \, Z_A, \, \bar{Z}_B  \, \rangle = \, \delta_{AB}, \quad \langle \,Z_0,
\, \bar{Z}_{N+1} \, \rangle = - \langle \,Z_{N+1}, \, \bar{Z}_0 \, \rangle=\, \mbox{i},
\end{equation}
while all other scalar products are zero. We define $\, p(Z) \, = \,
[Z_0]$, where $\, [Z_0]$ is the equivalence class of null vectors
represented by   $\, Z_0$. The left invariant Maurer-Cartan form
$\, \pi$ of SU$(N+1,1)$ is defined by the equation
\begin{equation}
d\,Z = Z \, \pi,\n
\end{equation}
which is in coordinates
\begin{equation}\label{struct1}
d (\, Z_0, \, Z_A, \, Z_{N+1} \, ) = ( \, Z_0, \, Z_B, \, Z_{N+1} \, ) \,
\begin{pmatrix}
\pi_0^0 & \pi_A^0 & \pi_{N+1}^0 \\
\pi_0^B & \pi_A^B & \pi_{N+1}^B \\
\pi_0^{N+1} & \pi_A^{N+1} & \pi_{N+1}^{N+1}
\end{pmatrix}.
\end{equation}
Coefficients of $\, \pi$ are subject to the relations obtained from
differentiating (\ref{prod}) which are
\begin{align}
\pi_0^0 + \bar{\pi}_{N+1}^{N+1} &=0 \n \\
\pi_0^{N+1}&= \bar{\pi}_0^{N+1}, \quad \pi_{N+1}^{0}= \bar{\pi}_{N+1}^{0} \notag \\
\pi_A^{N+1}&= - \, \mbox{i} \, \bar{\pi}_0^A, \quad \pi_{N+1}^A = \mbox{i} \, \bar{\pi}_A^0 \notag \\
\pi_B^A + \bar{\pi}_A^B &=0 \notag \\
\mbox{tr} \,  \pi &=0,  \n
\end{align}
and $\, \pi$ satisfies the structure equation
\begin{equation}\label{struct2}
- \, d\pi = \pi \wedge \pi.
\end{equation}

It is well known that the SU$(N+1,1)$-invariant CR structure on $\,
\Sigma^N \subset \mathbb{C}P^{N+1}$ as a real hypersurface is
biholomorphically equivalent to the standard  CR structure on
$\, S^{2N+1} = \partial \B^{N+1}$, where
$\, \B^{N+1} \subset \mathbb{C}^{N+1}$ is the unit ball.
The structure equation (\ref{struct1}) shows that for any local section
$\, s: \, \Sigma^N \, \to$ SU$(N+1,1)$, this CR structure is defined by
the hyperplane fields $\, (s^*\pi_0^{N+1})^{\perp} = \mathcal{H}$
and the set of (1,0)-forms $\, \{ \, s^*\pi_0^A \, \}$.

\vsp{1pc}
\tb{Definition.}
Let $\, M$ be a manifold of dimension $\, 2n+1$. A submanifold
defined  by an immersion $\, f: \, M \, \hookrightarrow \, \Sigma^N$
is a \emph{CR submanifold} if $\, f_*T_pM \, \cap \,
\mathcal{H}_{f(p)}$ is a complex subspace of $\, \mathcal{H}_{f(p)}$
of dimension $\, n$ for each $\, p \in M$.
\vsp{1pc}

\noindent Note the induced hyperplane fields $\, f_*^{-1}(f_*(TM) \, \cap \,
\mathcal{H})$ is necessarily a contact structure on $\, M$, and thus
$\, M$ has  an induced nondegenerate CR structure of hypersurface type.
When $\, M$ is equipped with a CR structure, an immersion
$\, f: M \to \Sigma^N$ is CR when the CR structure induced by $\, f$
is equivalent to the given one.

Consider $\, f^*$SU$(N+1,1) \to M$. From the definition, we may arrange so that
$\, \pi^{\alpha}_0=0$ for $\, n+1 \leq \alpha \leq N=n+m$ on this bundle.
Differentiating this,  we get
\begin{equation}
\pi^{\alpha}_i \w \pi^i_0 + \pi^{\alpha}_{N+1} \w \pi^{N+1}_0 = 0. \n
\end{equation}
By Cartan's lemma,
\begin{equation}\label{struct}
\pi^{\alpha}_i \equiv H^{\alpha}_{ij} \, \pi^j_0 \quad \mod \; \pi^{N+1}_0,
\end{equation}
for a coefficient $\, H^{\alpha}_{ij} = H^{\alpha}_{ji}$.
$\, H^{\alpha}_{ij}$ represents
the second fundamental form of $\, f$ [EHZ].

Our proof of \textbf{Theorem} is based on the following algebraic
theorem due to Iwatani on the asymptotic subspace of
the second fundamental form of a Bochner-K\"ahler submanifold [Iw][Br].
Let $\, V = \mathbb{C}^n$, $\, W = \mathbb{C}^{m}$
with the standard Hermitian scalar product.
Let $\,\{ \,  z^i \, \}$, $\, 1 \leq i \leq n$,
be a unitary (1,0)-basis for $\, V^*$,
and $\,\{ \,  w_{\ap} \, \} $, $\, 1 \leq \ap \leq m$,
be a unitary basis for $\, W$.
Let $\, S^{p,q}$ denote the space of polynomials of type $(p,q)$ on $\, V$.

\vsp{1pc}
\tb{Theorem [Iw]}\,
\emph{
Suppose $\, H = \, H^{\ap}_{ij} \, z^i z^j \otimes w_a
\in S^{2,0}\otimes W$ satisfies}
\begin{equation}
\gamma(H, H) = \, H^{\ap}_{ij} \, \bar{H}^{\ap}_{kl}
\, z^i z^j \otimes \bar{z}^k \bar{z}^l
= (z^k  \bar{z}^k ) \, h \in S^{2,2}, \quad h \in S^{1,1},\n
\end{equation}
\emph{or simply $\, \gamma(H, H)$ is Bochner-flat [Br].
Then the asymptotic vectors $\, \{ \, v \in V \; | \; H(v, v) = 0 \, \}$
form a subspace of $\, V$.
Let $\, n - \,\textnormal{k}$ be the dimension of this asymptotic subspace.
Then,
\begin{equation}
 m \geq \frac{1}{2}\, \textnormal{k} \,(2n-\,\textnormal{k}+1). \n
\end{equation}}

\noindent
We define k to be the \emph{rank} of $\, H$.
Note when $\, m < n$, \, k $=0$, and when $\, m < 2n-1$, \, k is at most 1.

In the case of our interest, the codimension of the CR submanifold is
bounded by $\, N-n = (3n-2)-n < 2n-1$, and hence k $\leq 1$.
Up to a unitary transformation on $\, V$, we may thus arrange
\begin{equation}
H^{\ap}_{ij}=\, H^{\ap}_{i} \, \delta_{jn}
 + H^{\ap}_{j} \, \delta_{in}, \n
\end{equation}
for coefficients $\, H^{\ap}_{i}$.
Set $\, \nu_i = H^{\ap}_{i} \, w_{\ap} \in W$.
A computation shows $\,\gamma(H, H)$ is Bochner-flat when
$
\langle\, \nu_i,  \, \nu_j \, \rangle = 0$ for  $\, i \ne j$, and
$\, \langle\, \nu_i,  \, \nu_i \, \rangle =  \, \langle\, \nu_j,  \, \nu_j \,\rangle.
$ for all $\, i, \, j$.
Up to a unitary transformation of $\, W$, we may set
\begin{align}
\nu_i &= \lambda \, w_i, \notag
\end{align}
for some $\, \lambda \geq 0$.

Let $\, f: \Sigma^n \hookrightarrow \Sigma^{N}$, $\, N=3n-2$,
be a local CR immersion.
Since $\, \Sigma^n$ is CR flat, after identifying
$\, V = f_*T^{\C}\Sigma^n \simeq \C^n \subset \mathcal{H}$ and
$\, W = V^{\perp} \simeq \C^{N-n} \subset \mathcal{H}$,
the second fundamental form of $\, f$ is Bochner-flat [EHZ].

Suppose $\, N-n < n$.
Then $\, H \equiv 0$, and $\, f$ is easily seen to be a part of the linear
embedding $\, f: \Sigma^n \hookrightarrow \Sigma^n \subset \Sigma^{N}$

Suppose $\, n \leq N-n < 2n-1$.
Then $\, H$ has rank at most 1, and from the argument above
we have in \eqref{struct},
\begin{align}
\pi^{n+i}_q &\equiv \lambda \, \delta_{iq} \, \pi_0^n
\quad \mod \; \pi^{N+1}_0, \quad \mbox{for} \,\; q < n, \; i \leq n, \label{seed} \\
\pi^{n+i}_n &\equiv \lambda (1+\delta_{in}) \,\pi_0^i
\quad \mod \; \pi^{N+1}_0, \quad \mbox{for} \,\; i \leq n, \n \\
\pi^{a}_i &\equiv 0
\quad \mod \; \pi^{N+1}_0, \quad \mbox{for} \; a > 2n, \n
\end{align}
for a coefficient $\, \lambda \geq 0$.

\section{Proof of theorem}
We wish to explore the consequence of the successive derivatives
of the relation \eqref{seed}.
To facilitate the computation, we shall agree on the index range
\begin{align}
1 \leq p, \, q, \, s, \, t \leq n-1,
\quad \; & p'=n + p, \n \\
1 \leq i,\, j,\, k,\, l \leq n,
\, \quad \quad \quad & i'=n + i,  \n\\
n'+1 \leq a, \, b \leq m=n'+r. \quad & \n
\end{align}
The induced Maurer-Cartan form \eqref{struct1} on $\, f^*$SU($N$+1,1)
decomposes according to these indices as follows.

\begin{equation}\label{decomp}
\pi =
\begin{pmatrix}
\pi_0^0 & \pi_q^0 &    \pi_n^0 &    \pi_{q'}^0 &       \pi_{n'}^0&   \pi_b^0 & \pi_{N+1}^0 \\
\pi_0^p & \pi_q^p &    \pi_n^p&     \pi_{q'}^p &        \pi_{n'}^p & \pi_b^p & \pi_{N+1}^p \\
\pi_0^n & \pi_q^n &    \pi_n^n &    \pi_{q'}^n &       \pi_{n'}^n &  \pi_b^n & \pi_{N+1}^n \\
\cdot &   \pi_q^{p'} & \pi_n^{p'} & \pi_{q'}^{p'} & \pi_{n'}^{p'} &  \pi_b^{p'} & \pi_{N+1}^{p'} \\
\cdot &   \pi_q^{n'} & \pi_n^{n'} & \pi_{q'}^{n'} & \pi_{n'}^{n'}&   \pi_b^{n'} & \pi_{N+1}^{n'} \\
\cdot &   \pi_q^a &    \pi_n^a &    \pi_{q'}^a &       \pi_{n'}^a &  \pi_b^a & \pi_{N+1}^a \\
\pi_0^{N+1} &  \pi_q^{N+1}& \pi_n^{N+1} & \cdot & \cdot   & \cdot                 & \pi_{N+1}^{N+1}
\end{pmatrix},
\end{equation}
where '$\cdot$' denotes $\, 0$.
We denote $\, \pi_0^i = \eta^i$, $\, \pi_0^{N+1} = \theta$, and
$\, -d \theta \equiv \, \mbox{i} \,
\eta^k \w \, \bar{\eta}^k \n = \, \mbox{i} \, \varpi \mod \theta$
\, for the sake of notation.
Equation \eqref{seed} for example can now be written as
\begin{equation}
\begin{pmatrix}
\pi_q^{p'} & \pi_n^{p'} \\
\pi_q^{n'} & \pi_n^{n'}  \\
\pi_q^a &    \pi_n^a
\end{pmatrix} \equiv
\begin{pmatrix}
\, \lambda \delta_{pq} \, \eta^n & \lambda \, \eta^p\\
\cdot      & 2 \lambda \, \eta^n  \\
\cdot      &   \cdot
\end{pmatrix} \mod{\theta}.\n
\end{equation}

The case $\, \lambda \equiv 0$($\, H \equiv 0$)
has already been identified with a part of the linear embedding.
Suppose $\, \lambda \ne 0$, and $\, H$ has rank 1.
We may scale $\, \lambda = 1$ using the group action by Re$\, \pi_0^0$,
and obtain the following normalized structure equation for
a nonlinear local CR immersion $\, f: \Sigma^n \hookrightarrow \Sigma^{N}$
with second fundamental form of rank 1,
\begin{equation}\label{rs1}
\begin{pmatrix}
\pi_q^{p'} & \pi_n^{p'} \\
\pi_q^{n'} & \pi_n^{n'}  \\
\pi_q^a &    \pi_n^a
\end{pmatrix} =
\begin{pmatrix}
 \delta_{pq} \, \eta^n &  \eta^p\\
\cdot      & 2  \eta^n  \\
\cdot      &   \cdot
\end{pmatrix} +
\begin{pmatrix}
h_q^{p'} & h_n^{p'} \\
h_q^{n'} & h_n^{n'}  \\
h_q^a    & h_n^a
\end{pmatrix} \theta,
\end{equation}
for coefficients $\, h^{i'}_j, \, h^a_j$.

\textbf{Theorem} is obtained by successive application of
Maurer-Cartan equation \eqref{struct2} to this structure equation.
We assume $\, n \geq 4$ for simplicity for the rest of this section,
as $\, n=3$ case can be treated with a minor modification.
The expression "differentiate X mod Y" would mean
"differentiate X and considering mod Y".

\vsp{1pc}
\emph{Step 1. $\;$}
Differentiate $\, \pi^{n'}_s  =  h^{n'}_s \, \theta$ $\mod{\theta}$, we get
\begin{equation}
\mbox{i}\,h^n_s \, \varpi \equiv ( \pi^{n'}_{s'} - 2 \pi^n_s) \w \eta^n
+ \pi^{n'}_{N+1} \w (-\, \im \, \bar{\eta}^s) \quad \mod{\theta}.  \n
\end{equation}
Since $\, n-1 \geq 2$, this implies $\, h^{n'}_s = 0$, and by Cartan's lemma
\begin{equation}
\begin{pmatrix}
\pi^{n'}_{s'} - 2 \pi^n_s \\
\pi^{n'}_{N+1}
\end{pmatrix}
\equiv
\begin{pmatrix}
2c_s & u \\
u  & 0
\end{pmatrix}
\begin{pmatrix}
\eta^n \\
-\im \, \bar{\eta}^s
\end{pmatrix}  \quad \mod{\theta}\n
\end{equation}
for coefficients $\, c_s, \, u$.

Differentiate $\, \pi^{a}_s  =  h^{a}_s \, \theta$ $\mod{\theta}$, we get
\begin{equation}
\mbox{i}\,h^a_s \, \varpi \equiv  \pi^{a}_{s'}  \w \eta^n
+ \pi^{a}_{N+1} \w (-\, \im \, \bar{\eta}^s) \quad \mod{\theta}.  \n
\end{equation}
Since $\, n-1 \geq 2$, this implies $\, h^{a}_s = 0$, and by Cartan's lemma
\begin{equation}
\begin{pmatrix}
\pi^{a}_{s'}  \\
\pi^{a}_{N+1}
\end{pmatrix}
\equiv
\begin{pmatrix}
2C^a_s & u^a \\
u^a  & 0
\end{pmatrix}
\begin{pmatrix}
\eta^n \\
-\im \, \bar{\eta}^s
\end{pmatrix}  \quad \mod{\theta}\n
\end{equation}
for coefficients $\, C^a_s, \, u^a$.

Differentiate $\, \pi^{a}_n  =  h^{a}_n \, \theta$ $\mod{\theta}$, we get
\begin{equation}
\mbox{i}\,h^a_n \, \varpi \equiv  \pi^{a}_{p'}  \w \eta^p + \pi^{a}_{n'}  \w 2 \eta^n
+ \pi^{a}_{N+1} \w (-\, \im \, \bar{\eta}^n) \quad \mod{\theta}.  \n
\end{equation}
Thus $\, h^a_n = u^a$, and
\begin{equation}
\pi^a_{n'} \equiv C^a_p \, \eta^p + C^a_n \eta^n - \im \, u^a \bar{\eta}^n
\mod{\theta} \n
\end{equation}
for coefficients $\, C^a_n$.

\emph{Step 2. $\;$}
Differentiate  $\, \pi^{t'}_s = h^{t'}_s \, \theta$ $\mod{\theta}\, $
for $\, t \ne s$, we get
\begin{equation}
\mbox{i}\,h^{t'}_s \, \varpi \equiv ( \pi^{t'}_{s'} -  \pi^t_s) \w \eta^n
- \pi^n_s \w \eta^t
+ \pi^{t'}_{N+1} \w (-\im \, \bar{\eta}^s)  \mod{\theta}.\n
\end{equation}
For $\, n-1 \geq 3$, this implies $\, h^{t'}_s = 0$ for $\, t \ne s$,
and by Cartan's lemma
\begin{equation}
\begin{pmatrix}
\pi^{t'}_{s'} -  \pi^t_s \\
-\pi^n_s \\
\pi^{t'}_{N+1}
\end{pmatrix}
\equiv
\begin{pmatrix}
0 & b_s   & -\im \, \bar{b}_t \\
b_s  & 0  & e \\
-\im \, \bar{b}_t  & e  & 0
\end{pmatrix}
\begin{pmatrix}
\eta^n \\
\eta^t \\
-\im \, \bar{\eta}^s
\end{pmatrix}  \mod{\theta}\n
\end{equation}
for coefficients $\, b_s, \, e$.
Since $\, \pi^{t'}_{s'} -  \pi^t_s$ is skew Hermitian, it cannot have
any $\, \eta^n$-term.

\emph{Step 3. $\;$}
Differentiate  $\, \pi^{t'}_t = \eta^n + h^{t'}_t \, \theta \,$
$\mod{\theta}$, we get
\begin{equation}
h^{t'}_t \, \varpi \equiv
\Delta_t \w \eta^n
+( b_p \, \eta^n - \im \, e  \, \bar{\eta}^p ) \w \eta^p
+( -b_t \, \eta^t + \bar{b}_t \, \bar{\eta}^t ) \, \w \eta^n  \mod{\theta}, \n
\end{equation}
where $\, \Delta_t = \pi^{t'}_{t'} -  \pi^t_t+ \pi^0_0 - \pi^n_n$.
This implies $\, h^{t'}_t = e$, and
\begin{equation}
\Delta_t =  a_t \, \eta^n - \im \, e \, \bar{\eta}^n
+ ( b_t \, \eta^t - \bar{b}_t \, \bar{\eta}^t )
+ \sum_p \, b_p \,  \eta^p  - A_t \, \theta, \n
\end{equation}
for coefficients $\, a_t, \, A_t$.

\emph{Step 4. $\;$}
Since $\, h^{t'}_s = \delta_{st} \,e$ and from \eqref{rs1},
we may add $\, e \, \theta$  to
$\, \eta^n$ to translate $\, e = 0$, which we assume from no on.
We also translate $\, h^{t'}_n = 0\,$ similarly by adding
$\, h^{t'}_n \, \theta$ to $\, \eta^t$.
Differentiating
$\, \pi^{t'}_n = \eta^t$ $\mod{\theta}\, $ with these relations
and collecting terms, we get $\, b_t = c_t = 0$, and
\begin{equation}
0 \equiv a_t \,\eta^n \w \eta^t - 2 \im \, \bar{u} \,\eta^t \w \eta^n
\mod{\theta}. \n
\end{equation}
Thus $\, a_t = - 2 \,\im \,\bar{u}$.

\emph{Step 5.} $\;$ Differentiate
$\, \pi^{n'}_n = 2 \, \eta^n + h^{n'}_n \, \theta$ $\mod{\theta}$,
and collecting terms, we get $\, h^{n'}_n = u$, and
$
0 \equiv \Delta_n \w \eta^n - \im u \, \eta^n \w \bar{\eta}^n \mod{\theta},
$
where $\, \Delta_n  = \pi^{n'}_{n'} -  \pi^n_n+ \pi^0_0 - \pi^n_n $.
Hence
\begin{equation}
\Delta_n = -\im u\, \bar{\eta}^n + a_n \, \eta^n -A_n \, \theta, \n
\end{equation}
for coefficients  $\, a_n, \, A_n$.
But $\, \Delta_t-\Delta_n$ is purely imaginary, and comparing with
\emph{Step 3}, $\, a_n = - 3 \,\im \,\bar{u}$.

\emph{Step 6. $\;$}   Now by considering $\, \theta$-terms in
\emph{Step 1, 2, 3, 4, 5}, and the fact that
$\, \pi^{\ap}_i \w \eta^i + \pi^{\ap}_{N+1} \w \theta = 0$,
we obtain the following simple structure equations.
We omit the details of computations.

\begin{align}
\begin{pmatrix}
\pi^{p'}_q  & \pi^{p'}_n \\
\pi^{n'}_q  & \pi^{n'}_n \\
\pi^{a}_q  &  \pi^{a}_n \\
\end{pmatrix}
&=
\begin{pmatrix}
\delta_{pq} \, \eta^n  &  \eta^p \\
0  & 2 \, \eta^n + u \, \theta\\
0  &  u^a \, \theta
\end{pmatrix}\n \\
\begin{pmatrix}
\pi^{p'}_{N+1}   \\
\pi^{n'}_{N+1}   \\
\pi^{a}_{N+1}    \\
\end{pmatrix}
&=
\begin{pmatrix}
0\\
u \, \eta^n\\
u^a \, \eta^n
\end{pmatrix}\n \\
\begin{pmatrix}
\pi^{n'}_{s'} \\
\pi^{a}_{s'} \\
\end{pmatrix}
&=
\begin{pmatrix}
-\im u\, \bar{\eta}^s   \\
-\im u^a\, \bar{\eta}^s + 2C^a_s \, \eta^n
\end{pmatrix}\n \\
\pi^n_s &= 0,\, \pi^{t'}_{s'} = \pi^t_s \; \; \mbox{for} \; t \ne s,\n \\
\pi^a_{n'} &= C^a_p \, \eta^p + C^a_n \, \eta^n -\im u^a \, \bar{\eta}^n
+ h^a_{n'} \, \theta, \n\\
\begin{pmatrix}
\pi^t_{N+1} \\
\pi^n_{N+1}
\end{pmatrix} &=
\begin{pmatrix}
(A - \im \, ( u \, \bar{u} + u^a \, \bar{u}^a )) \, \eta^t
-2u^a \bar{C}^a_t \, \bar{\eta}^n + B_t \, \theta   \\
A \, \eta^n + B_n \, \theta
\end{pmatrix} \n \\
\begin{pmatrix}
\Delta_t \\
\Delta_n
\end{pmatrix} &=
\begin{pmatrix}
- 2 \,  \im  \,\bar{u} \, \eta^n  - A \, \theta \\
 -3 \,  \im  \, \bar{u}  \, \eta^n -  \,  \im  \, u \, \bar{\eta^n} \,
- A_n \, \theta
\end{pmatrix} \n \\
A_n + \bar{A}_n &= A + \bar{A} \n \\
du &= u \, ( \pi^{N+1}_{N+1}-\pi^0_0 +\pi^n_n-\pi^{n'}_{n'})
+2(A-A_n)\, \eta^n - u^a \, \pi^{n'}_a + u_0 \, \theta \n \\
du^a &= u^a \, ( \pi^{N+1}_{N+1}-\pi^0_0 +\pi^n_n) - u^b \, \pi^a_b
- u \, \pi_{n'}^a + 2 h^a_{n'} \, \eta^n  + u^a_0 \, \theta  \label{ua}
\end{align}

\emph{Proof of Theorem}.$\;$
From [HJX, Corollary 1.4], a proper holomorphic map
$\, F: \B^{n+1} \to \B^{3n-1}$ which is $\, C^3$ up to the boundary
is rational, in particular real analytic up to the boundary.
Hamada recently showed that a rational proper holomorphic map
from $\, \B^{n+1}$ to $\, \B^{2n+2}$, $\, n \geq 3$, is equivalent to
one of the three class of maps in   \tb{Theorem}.
It thus suffices to show that any analytic local CR immersion
$\, f: \Sigma^{n} \hook  \Sigma^{3n-2}$ in fact lies in
$\, f:  \Sigma^{n} \hook  \Sigma^{2n+1} \subset  \Sigma^{3n-2}$
for $\, n \geq 3$.

\emph{Step 1. $\;$}
From the refined structure equations above, differentiate
$\, \pi^{p'}_{N+1} =0$, $\, \pi^{n'}_{N+1} = u \, \eta^n$,
and $\, \pi^{a}_{N+1} = u^a \, \eta^n$. After a short computation,
we get  $\, B_p =0, \, B_n = 0$, and
\begin{equation}\label{971}
u_0 = -2 u A, \; u^a_0 = -2 u^a A.
\end{equation}
Differentiating $\, \pi^{n}_{q} = 0$, and $\, \pi^{t'}_{s'} = \pi^t_s$
for $\, t \ne s$ with these relations, we get
\begin{align}
A - \bar{A} &= \im \, ( u \, \bar{u} + u^a \, \bar{u}^a - 1 ), \n \\
\sum_{a=n'+1}^{n'+r} C^a_t \bar{C}^a_s &= 0, \quad \mbox{for} \; \; t \ne s.
\label{972}
\end{align}

\emph{Step 2. $\;$}
Differentiate $\, \im \, \pi^{n'}_{q'} = u \, \bar{\eta}^q$,
and collecting $\, \eta^n \w \bar{\eta}^q$-terms, we get
\begin{equation}
\sum_{a=n'+1}^{n'+r} C^a_q \bar{C}^a_q = 1- u \bar{u}+\im (A-A_n), \n
\end{equation}
which is independent of the index $\, q$.
Since $\, r < n-1$ from our assumption on the codimension,
this and \eqref{972} force
\begin{equation}\label{975}
C^a_q = 0,
\end{equation}
and $\, A_n = A + \im ( u \bar{u} -1)$.

\emph{Step 3. $\;$}
Differentiate
$\, \pi^a_{n'} = C^a_n \, \eta^n -\im u^a \, \bar{\eta}^n + h^a_{n'} \, \theta$
$\mod{ \theta, \, \eta^n, \ \bar{\eta}^n}$, we get
\begin{equation}\label{973}
h^a_{n'} = -\im u^a  \bar{u}.
\end{equation}

\emph{Step 4. $\;$}
Differentiate
$\, \pi^p_{N+1} = (\, \bar{A} - \im) \, \eta^p$,
$\, \pi^n_{N+1} = A \, \eta^n$, we get
\begin{equation}
dA=  A \, (\, \pi^{N+1}_{N+1} - \pi^0_0) +\pi^0_{N+1}
+2(u\, \bar{\eta}^n - \bar{u} \, \eta^n) +
(u \bar{u} + u^a \bar{u}^a-A^2) \, \theta.\n
\end{equation}
Note
$\,
\Delta_t + \bar{\Delta}_t  = \pi^0_0 + \bar{\pi}^0_0 = \pi^0_0 - \pi^{N+1}_{N+1}
 =2 \,\im \, ( u\, \bar{\eta}^n - \bar{u} \, \eta^n) - (A + \bar{A}) \, \theta$.
Now finally differentiating
$\, \Delta_n = -3 \,  \im  \, \bar{u}  \, \eta^n -  \,  \im  \, u \, \bar{\eta^n} \,
- A_n \, \theta$ with these relations,
and collecting $\, \eta^n \w \bar{\eta}^n$-terms,
we get $\, \sum_{a=n'+1}^{n'+r}  C^a_n \bar{C}^a_n=0$,
and hence
\begin{equation}
C^a_n =0. \label{974}
\end{equation}

Case \tb{B}.\;
Suppose $\, u^a =0$ for all $\, a$. From \eqref{975}, \eqref{973}, \eqref{974},
$\, \pi^a_n = \pi^a_{s'} =\pi^a_{n'} =\pi^a_{N+1} = 0$.
In the notation of \eqref{struct1},
the complex $\, (2n+2)-$plane $\, Z_0 \w Z_1 \w ... \, Z_{2n} \w Z_{N+1}$
is then constant along the CR immersion $\, f$.
Hence $\, f: \Sigma^n \hook \Sigma^{2n} \subset \Sigma^{3n-2}$.

Case \tb{C}.\;
Suppose $\,\vec{u} = (u^{n'+1}, \, ... \, u^{n'+r}) \ne 0$.
Using a group action by $\, \pi^a_b$, we may rotate so that
$\,\vec{u} = (u^{n'+1}, \, 0, \, ... \, 0)$ with $\, u^{n'+1} \ne 0$.
But from \eqref{975}, \eqref{973}, \eqref{974}, we have
$\, \pi^a_n = \pi^a_{s'} =\pi^a_{n'} =\pi^a_{N+1} = 0$ for $\, a > n'+1$.
Moreover, \eqref{ua} shows $\, \pi_{n'+1}^{a} = 0$ for $\, a > n'+1$.
In the notation of \eqref{struct1}, the complex $\, (2n+3)-$plane
$\, Z_0 \w Z_1 \w ... \, Z_{2n} \w Z_{2n+1} \w Z_{N+1}$
is then constant along the CR immersion $\, f$.
Hence $\, f: \Sigma^n \hook \Sigma^{2n+1} \subset \Sigma^{3n-2}$.\sq

At this stage, note that the only possibly independent coefficients in the
structure equations are $\, A, \, u, \, u^a$, and that the expression
for their derivatives does not involve any new variables.
The structure equations for local CR immersion
$\, f: \Sigma^n \hook \Sigma^{3n-2}$ thus close up at order 3.
A long but direct computation shows that these equations are compatible,
i.e.,   $\, d^2 = 0$ is a formal identity of the structure equation.
In [Wa], we showed that every local CR immersion
$\, f: \Sigma^n \hook \Sigma^{2n}$
is equivalent to a part of either the linear embedding or Whitney map.
We suspect a similar argument can be applied to show that
every linearly full local CR immersion $\, f: \Sigma^n \hook \Sigma^{2n+1}$
is equivalent to a boundary CR map of a type \tb{C} proper holomorphic map
in \tb{Theorem}.

\Remark The computation involved here is reminiscent of Cartan's
local isometric embedding of Hyperbolic space $\, \mathbb{H}^n$ in
Euclidean space $\, \mathbb{E}^{2n-1}$ via exteriorly orthogonal
symmetric bilinear forms [Ca].
Overdetermined nature of CR geometry forces the structure equation
to close up instead of becoming involutive.

\Remark It is natural to ask the K\"ahler analogue of this rigidity theorem.
Let $\, M^n \hook X^N_{\epsilon}$ be a complex submanifold in a complex
space form $\, X$ of constant holomorphic sectional curvature $\, \epsilon$.
Suppose the induced metric on $\, M$ is Bochner-K\"ahler and
the codimension is bounded by $\, N - n < 2n-1$.
A short computation shows that such $\, M$ is totally geodesic for $\, n \geq 3$.

The computation suggests that the gap phenomena in \tb{Theorem}
may persist for CR immersions between spheres
with second fundamental form of rank 1;   \emph{
linearly full CR immersions $\, f: \Sigma^n \hook \Sigma^{N}$,
$\, N \leq n^2-2$, with second fundamental form of rank 1 can only occur
when $\, N = \mu   n, \, \mu  n +1, \, ... \,, \,  \mu   n + \mu-1$
for $\, \mu \leq n-1$}.

Consider the following generalization of Whitney map.
\begin{align}
F&: \, \B^{n+1} \to \B^{\mu \, (n + 1)}, \label{list} \\
F&(z^i, \, z^0) =
(x_1 \, z^i, \, x_2 \, z^i z^0, \, ... \, ,\, x_{\mu} \, z^i (z^0)^{\mu-1}, \,
y_1 \, z^0, \, y_2 \, z^0 z^0, \, ... \, ,\, y_{\mu} \, z^0 (z^0)^{\mu-1}). \n
\end{align}
where $\, x_A, \, y_A$ are constants satisfying
$\, x_1=1, \, x_{\mu} = y_{\mu}, \, x^2_A = y^2_A + x^2_{A+1}$
for $\, 1 \leq A \leq \mu-1$.
The geometric rank of a proper holomorphic map between unit balls
introduced by Huang is likely bounded above by the rank of
the second fundamental form of its boundary CR map at a generic point [Hu].
It would be interesting to understand how exhaustive \eqref{list} is
among the set of proper holomorphic maps with geometric rank 1
from $\, \B^{n+1}$  to $\,\B^{\mu \, (n + 1)}$  for $\, \mu \leq n-1$.

\vsp{2pc}
\noindent
\begin{center}
\textbf{\large{References}}
\end{center}

\noindent [Br] Bryant, Robert, \emph{
Bochner-K\"ahler metrics},
J. Amer. Math. Soc.  14  (2001)   no. 3, 623--715

\noindent [Ca] Cartan, E., \emph{
Sur les varietes de courbure constante
d'un espace euclidean ou non euclidean},
Bull. Soc. Math. France, 47 (1919),  48 (1920)  125--160, 32--208

\noindent [ChM]  Chern, S. S.; Moser, J. K., \emph{
Real hypersurfaces in complex manifolds},
S. S. Chern Selected papers vol III, 209--262

\noindent [DA1]  D'Angelo, John P., \emph{
Proper holomorphic maps between balls of different dimensions},
Michigan Math. J.  35 (1988)  83--90

\noindent [DA2]  \underline{\quad \quad}, \emph{
Several Complex Variables and the Geometry of Real Hypersurfaces},
Studies in Advanced Mathematics, CRC Press, 1993

\noindent [EHZ]  Ebenfelt, P.; Huang, Xiaojun; Zaitsev, Dmitri, \emph{
Rigidity of CR-immersions into Spheres},
Comm. Anal. Geom.  12 (2004)  no. 3, (2004) 631--670

\noindent [Fo]  Forstneri\v c, F.,\emph{
Proper holomorphic maps between balls},
Duke Math. J. 53 (1986)  427--441

\noindent [Ha]  Hamada, H., \emph{
Rational proper holomorphic maps from
$\,\bold B\sp n$ into $\,\bold B\sp {2n}$},
Math. Ann.  331  (2005)  693--711

\noindent [Han]  Han, Chong-Kyu, \emph{
Complete differential system for the mappings of CR manifolds
of nondegenerate Levi forms},
Math. Ann.  309  (1997)   no. 3, 401--409

\noindent [HJ]   Huang, Xiaojun; Ji, Shanyu, \emph{
Mapping $\,\bold B\sp n$ into $\,\bold B\sp {2n-1}$},
Invent. Math.  145  (2001)   no. 2, 219--250

\noindent [HJX]  \underline{\quad \quad}, \underline{\quad \quad},
Xu, D., \emph{
Several Results for Holomorphic Mappings from
$\,\bold B\sp n$ into $\,\bold B\sp {N}$},
Geometric analysis of PDE and several complex variables,
Contemp. Math., 368, AMS, (2005) 267--292

\noindent [Hu] Huang, Xiaojun,  \emph{
On a semi-rigidity property for holomorphic maps},
Asian J. Math.  7 (2003) 463--492

\noindent [Iw]  Iwatani, Teruo, \emph{
K\"ahler submanifolds with vanishing Bochner curvature tensor},
Mem. Fac. Sci. Kyushu Univ. Ser. A 30 (1976)  no. 2, 319--321

\noindent [Wa]  Wang, S. H., \emph{
1-rigidity for CR submanifolds in spheres},
arXiv. math.DG/0506134

\vspace{2pc}

\noindent

\noindent Sung Ho Wang \\
\noindent Department of Mathematics \\
\noindent Kias \\
\noindent Seoul, Corea 130-722 \\
\texttt{shw@kias.re.kr}

\end{document}